\input amssym.def
\input amssym
\documentstyle[12pt]{article}
\newcommand{\doublespace}{
   \renewcommand{\baselinestretch}{1.5}
   \large\normalsize}

\def \Z{\Bbb Z}
\def \C{\Bbb C}

\def \Q{\Bbb Q}
\def \N{\Bbb N}
\def \wt{{\rm wt}}

\def \Res{{\rm Res}}

\def \End{{\rm End}}

\def \Aut{{\rm Aut}}

\def \Vir{{\rm Vir}}
\def \<{\langle}
\def \>{\rangle}

\def \o{\omega }

\def \l{\lambda }

\doublespace
\begin{document}
\newtheorem{th1}{Theorem}
\newtheorem{th}{Theorem}[section]
\newtheorem{prop}[th]{Proposition}
\newtheorem{coro}[th]{Corollary}
\newtheorem{lem}[th]{Lemma}
\newtheorem{rem}[th]{Remark}
\newtheorem{de}[th]{Definition}

\begin{center}
{\Large {\bf Introduction to vertex operator algebras I}} \\
\vspace{0.5cm}
Chongying Dong\footnote{Supported by NSF grant
DMS-9303374 and a research grant from the Committee on Research,
UC Santa Cruz.} \\
Department of Mathematics, University of
California, Santa Cruz, CA 95064
\end{center}

\hspace{1.5 cm}

\section{Introduction}

The theory of vertex (operator) algebras has developed rapidly in the last few
years. These rich algebraic structures
provide the proper formulation for the moonshine module
construction for the Monster group ([B1-B2], [FLM1], [FLM3]) and also give a
lot of  new insight into the
representation theory of the Virasoro algebra and affine Kac-Moody algebras
(see for instance [DL3], [DMZ], [FZ], [W]).
The modern notion of chiral algebra in  conformal field theory [BPZ]
in physics essentially
corresponds to the mathematical
notion of vertex operator algebra; see e.g. [MS].

This is the first part of three consecutive lectures by Huang, Li and myself.
In this part we are mainly concerned with the definitions of vertex operator
algebras, twisted modules and examples. The second part by Li is about
the duality and local systems and the third part by Huang is devoted to
the contragradient modules and geometric interpretations of vertex operator
algebras. (We refer the reader to Li and Huang's lecture notes for the
related topics.) So
 many exciting topics are not covered in these three lectures.
The books [FLM3] and [FHL] are excellent introduction to the subject.
There are also existing papers [H1], [Ge] and [P] which review the axiomatic
definition of vertex operator algebras, geometric interpretation
of vertex operator algebras, the connection with conformal
field theory, Borcherds algebras and the monster Lie algebra.

Much work on vertex operator algebras has been concentrated on the concrete
examples of vertex operator algebras and the representation theory.
In particular, the representation theory for the concrete vertex operator
algebras, which include
the moonshine vertex operator algebra $V^{\natural}$ ([FLM3],[D3]),
the vertex operator algebras based on even positive definite lattices [D1],
the vertex operator algebras associated with the integrable representations
of affine Lie algebras, Virasoro algebra and $W$-algebras
([DMZ], [DL3], [FKRW], [FZ], [KW] and [W]), have
been studied extensively. There are also abstract approaches such as
one to one correspondence between the set of
inequivalent irreducible (twisted) modules for a
given vertex operator algebra and  the set of inequivalent
irreducible modules for an
associative algebra associated with the vertex operator algebra and
an automorphism of the algebra (see [DLM1] and [Z]),
the theory of local system [L1-L2],
induced modules [DLi], the tensor products of modules ([HL1-HL5], [H3]
 and [L3]), the classification of certain vertex operator algebras
[Lia]; See also [FHL]
for the results concerning intertwining operators and contragradient
modules.
Many of these results are analogues of the corresponding results in
the classical Lie algebra theory, usually much more subtle.

This paper is organized as follows: In Section 2, we present the definitions
of vertex operator algebras and twisted modules.
We also make some remarks.
The Section 3 is about the examples of vertex operator algebras.
In particular, we discuss the vertex operator algebras associated with
Heisenberg algebras, positive definite even lattices, affine Lie algebras
and Virasoro algebra. We emphasize how a vertex  operator algebras with
finitely many generators can be constructed by using the Jacobi identity.
In Section 5 we finally mention the orbifold theory which is not treated
in the introductory text and various generalizations of the notion of
vertex operator algebras.

{\bf Acknowledgment.} This paper is an expanded version of my lecture
in the workshop of ``Moonshine and vertex operator algebra'' at Research
Institute of Mathematical Science at Kyoto in the Fall of 1994. I
thank the organizer Professor Miyamoto for the opportunity to present talk at
this stimulating workshop. I also would like to thank Professors Miyamoto and
Yamada for their hospitality during my visit in Japan.

\section{Vertex operator algebras and modules}
\setcounter{equation}{0}
The notion of vertex (operator) algebra arose naturally from the problem
of realizing the monster simple group as a symmetry group of certain
algebraic structure. In this section we shall review the definition of vertex
operator algebras and
their (twisted) modules (see [B1], [D2], [FFR], [FHL] and [FLM3]).
We shall also discuss some consequences of the definitions and basic
properties.

First we introduce some notation. We shall use  commuting formal
variables $z,$ $z_0,$ $z_1,$ $z_2,$ etc., and the basic generating function
\begin{equation}
\delta(z)=\sum_{n\in{\Bbb Z}}z^n,
\end{equation}
formally the expansion of the $\delta$-function at $z=1.$ The
fundamental (and elementary)
properties of the $\delta$-function can be found in [FLM3], [FHL] and [DL3].

For a vector space $V$ (we work over ${\Bbb C}$) and a positive integer $r$
we denote the vector space of formal Laurent series
in $z^{1/r}$ with coefficients in $V$ by
\begin{equation}\label{oo}
V[[z^{1/r},z^{-1/r}]]=\left\{\sum_{n\in\frac{1}{r}\Z}v_nz^n|v_n\in V\right\}.
\end{equation}
We also define a formal residue notation:
$$\Res_z\sum_{n\in\frac{1}{r}\Z}v_nz^n=v_{-1}.$$

A {\em vertex  operator algebra} is a ${\Z}$-graded vector space:
\begin{equation}
V=\coprod_{n\in{\Z}}V_n; \ \ \ \mbox{for}\ \ \ v\in V_n,\ \
n=\mbox{wt}\,v;
\end{equation}
such that $\dim V_n<\infty$ for all $n\in\Z$ and $V_n=0$ if $n$ is sufficiently
small; equipped with a linear map
\begin{equation}\label{3.30}
\begin{array}{l}
V \to (\mbox{End}\,V)[[z,z^{-1}]]\\
v\mapsto Y(v,z)=\displaystyle{\sum_{n\in{\Z}}v_nz^{-n-1}}\ \ \ \  (v_n\in
\mbox{End}\,V)
\end{array}
\end{equation}
and with two distinguished vectors ${\bf 1}\in V_0,$ $\omega\in V_2$
satisfying the following conditions for $u, v \in V$:
\begin{eqnarray}
& &u_nv=0\ \ \ \ \ \mbox{for}\ \  n\ \ \mbox{sufficiently large};\\
& &Y({\bf 1},z)=1;\\
& &Y(v,z){\bf 1}\in V[[z]]\ \ \ \mbox{and}\ \ \ \lim_{z\to
0}Y(v,z){\bf 1}=v;
\end{eqnarray}
\begin{equation}\label{jac}
\begin{array}{c}
\displaystyle{z^{-1}_0\delta\left(\frac{z_1-z_2}{z_0}\right)
Y(u,z_1)Y(v,z_2)-z^{-1}_0\delta\left(\frac{z_2-z_1}{-z_0}\right)
Y(v,z_2)Y(u,z_1)}\\
\displaystyle{=z_2^{-1}\delta
\left(\frac{z_1-z_0}{z_2}\right)
Y(Y(u,z_0)v,z_2)}
\end{array}
\end{equation}
(Jacobi identity) where all binomial
expressions, for instance, $(z_1-z_2)^n\ (n\in {\Z})$ are to be expanded
in nonnegative
integral powers of second variable $z_2$: This identity is
interpreted algebraically as follows: if this identity is applied to a
single vector of $V$ then the coefficient of each monomial in
$z_0,z_1,z_2$ is a finite sum in $V;$
\begin{equation}
[L(m),L(n)]=(m-n)L(m+n)+\frac{1}{12}(m^3-m)\delta_{m+n,0}(\mbox{rank}\,V)
\end{equation}
for $m, n\in {\Z},$ where
\begin{equation}
L(n)=\omega_{n+1}\ \ \ \mbox{for}\ \ \ n\in{\Z}, \ \ \
\mbox{i.e.},\ \ \ Y(\omega,z)=\sum_{n\in{\Z}}L(n)z^{-n-2}
\end{equation}
and
\begin{eqnarray}
& &\mbox{rank}\,V\in {\Q};\\
& &L(0)v=nv=(\mbox{wt}\,v)v \ \ \ \mbox{for}\ \ \ v\in V_n\\
& &\mbox{rank}\,V\in {\Q};\\
& &L(0)v=nv=(\mbox{wt}\,v)v \ \ \ \mbox{for}\ \ \ v\in V_n\
(n\in{\Z}); \label{3.40}\\
& &\frac{d}{dz}Y(v,z)=Y(L(-1)v,z).\label{3.41}
\end{eqnarray}
This completes the definition. We denote the vertex
operator algebra just defined by
\begin{equation}
(V,Y,{\bf 1},\omega)
\end{equation}
(or briefly, by $V$). The series $Y(v,z)$ are called vertex operators.

The following are consequences of the definitions:
For $u,v\in V,$
\begin{eqnarray}
& &[L(-1),Y(v,z)]=Y(L(-1)v,z)\label{6.20}\\
& &[L(0),Y(v,z)]=Y(L(0)v,z)+zY(L(-1)v,z)\label{6.21}\\
& &Y(e^{z_0L(-1)}v,z)=Y(v,z+z_0)\label{6.26}\\
& &e^{z_0L(-1)}Y(v,z)e^{-z_0L(-1)}=Y(e^{z_0L(-1)}v,z)\label{6.27}\\
& &Y(v,z)\mbox{\bf 1}=e^{zL(-1)}v\label{6.28}\\
& &Y(u,z)v=e^{zL(-1)}Y(v,-z)u.\label{6.29}
\end{eqnarray}

Let $S$ be a subset of $V.$ The subalgebra $\<S\>$ {\em generated} by $S,$
defined as
the smallest subalgebra containing $S,$ is given by
$$\{v_{n_2}^1\cdots v_{n_k}^k\cdot {\bf 1} |v^i\in S\cup \{{\bf
1},\omega\},n_i\in \Z\}.$$  A vertex operator algebra $V$ is generated by $S$
if
$V=\<S\>.$ If we further assume that $S$ is a finite set we say that
$V$ is finitely generated.

An {\em ideal} $I$ of $V$ is a subspace of $V$ such that ${\bf 1}\not\in I$ and
$u_nv\in I$ for any $u\in V$ and $v\in I.$ Using the skew symmetry (\ref{6.29})
one can verify that the quotient space $V/I$ has a structure of vertex
operator algebra.

An {\em automorphism} $g$ of the vertex operator
algebra $V$ is a linear automorphism of $V$ preserving ${\bf 1}$ and $\omega$
such that the actions of $g$
and $Y(v,z)$ on $V$ are compatible in the sense that
\begin{equation}\label{aut}
gY(v,z)g^{-1}=Y(gv,z)
\end{equation}
for $v\in V.$ Then $gV_n\subset V_n$ for $n\in\Z$ and $V$ is a direct sum of
the eigenspaces of $g:$
\begin{equation}\label{dec}
V=\coprod_{j\in \Z/r\Z}V^j
\end{equation}
where $r$ is the order of $g,$ $\eta=e^{2\pi i/r}$ and $V^j=\{v\in
V|gv=\eta^jv\}.$

We also have the notion of $g$-twisted module
(see [D2], [FFR], [Le] and [FLM2]).
Let $(V,Y,{\bf 1},\omega)$ be a vertex
algebra and let $g$ be an automorphism of $V$ of order $r.$ A $g$-twisted
$module$ $M$ for $(Y,V,{\bf 1},\omega)$ is a $\C$-graded
vector space:
\begin{equation}\label{6.62}
M=\coprod_{n\in{\C}}M_n; \ \ \ \mbox{for}
\ \ \ w\in M_n,\ \ n=\mbox{wt}\,w;
\end{equation}
such that for any fixed $\lambda\in\C,$ $M_{\l+n}=0$ if
$n\in\frac{1}{r}\Z$ is sufficiently small and $\dim M_c<\infty$ for all
$c\in \C;$ equipped with a linear map
\begin{equation}
\begin{array}{l}
V\to (\mbox{End}\,M)[[z^{1/r},z^{-1/r}]]\label{map}\\
v\mapsto\displaystyle{ Y(v,z)=\sum_{n\in{\frac{1}{r}\Z}}v_nz^{-n-1}\ \ \
(v_n\in
\mbox{End}\,M)}
\end{array}
\end{equation}
satisfying the following conditions for $u,v\in V$ and
$w\in M$ and $l\in\Z$:
\begin{eqnarray}
& &Y(v,z)=\sum_{n\in l/r+\Z}v_nz^{-n-1}\ \ \ \ {\rm for}\ \ v\in
V^l;\label{1/2}\\
& &u_nw=0\ \ \ \mbox{for}\ \ \ n\ \ \ \mbox{sufficiently\ large};\\
& &Y({\bf 1},z)=1;
\end{eqnarray}
\begin{equation}\label{jacm}
\begin{array}{c}
\displaystyle{z^{-1}_0\delta\left(\frac{z_1-z_2}{z_0}\right)
Y(u,z_1)Y(v,z_2)-z^{-1}_0\delta\left(\frac{z_2-z_1}{-z_0}\right)
Y(v,z_2)Y(u,z_1)}\\
\displaystyle{=z_2^{-1}\left(\frac{z_1-z_0}{z_2}\right)^{-i/r}
\delta\left(\frac{z_1-z_0}{z_2}\right)Y(Y(u,z_0)v,z_2)}
\end{array}
\end{equation}
where $u\in V^i$ and $Y(u,z_0)$ is an operator on $V$;
\begin{equation}
[L(m),L(n)]=(m-n)L(m+n)+\frac{1}{12}(m^3-m)\delta_{m+n,0}(\mbox{rank}\,V)
\end{equation}
for $m, n\in {\Z},$ where
\begin{eqnarray}
& &L(n)=\omega_{n+1}\ \ \ \mbox{for}\ \ \ n\in{\Z}, \ \ \
\mbox{i.e.},\ \ \ Y(\omega,z)=\sum_{n\in{\Z}}L(n)z^{-n-2};\label{6.70}\\
& &L(0)w=nw=(\mbox{wt}\,w)w \ \ \ \mbox{for}\ \ \ w\in M_n\
(n\in{\Q});\label{6.71}\\
& &\frac{d}{dz}Y(v,z)=Y(L(-1)v,z).\label{6.72}
\end{eqnarray}
This completes the definition. We denote this module by
$(M,Y)$ (or briefly by $M$).

\begin{rem} If $g=1,$ then $M$ is an ordinary
module in the precise sense of [FLM3]. Note that any $g$-twisted $V$-module
is a $V^0$-module.
\end{rem}

\begin{rem} In the definition of twisted module if removing the
finiteness condition $\dim M_c<\infty$ we get a {\em weak} $g$-twisted module.
\end{rem}

\begin{rem} Taking $\Res_{z_0}$ in (\ref{jacm}), we get the commutator formula:
\begin{equation}\label{com}
[Y(u,z_1),Y(v,z_2)]=
\Res_{z_0}\left\{z_2^{-1}\left(\frac{z_1-z_0}{z_2}\right)^{-i/r}
\delta\left(\frac{z_1-z_0}{z_2}\right)Y(Y(u,z_0)v,z_2)\right\}.
\end{equation}
Note that the factor $\left(\frac{z_1-z_0}{z_2}\right)^{-i/r}
\delta\left(\frac{z_1-z_0}{z_2}\right)$ only involves the nonnegative
powers of $z_0.$ Thus in the computation of commutator $[Y(u,z_1),Y(v,z_2)]$
we only use the ``singular terms'' in $Y(u,z_0)v,$ namely,
$\sum_{n\geq 0}u_nvz_0^{-n-1}.$ This fact was well-known in the physics
literature. Moreover, for $ s, t  \in \Q $, we can compare the coefficients of
$ z_{1}^{-s-1}z_{2}^{-t-1}$ on the both sides of (\ref{com}) to get
the commutator $[u_s,v_n]:$
\begin{equation}
[u_{s},v_{t}]=u_sv_t-v_tu_s=
\sum_{m\geq 0}{s\choose m}(u_{m}v)_{s+t-m}.
\label{jacobi-term}
\end{equation}
\end{rem}

\begin{rem} Acting on $M$ the operator $Y(Y(u,z_0)v,z_2)$ is determined
uniquely by the operators $Y(u,z_1)$ and $Y(v,z_2).$
In order to see this, we first recall a relation on $\delta$-functions
{}from [FLM3]:
$$z_1^{-1}\left(\frac{z_2+z_0}{z_1}\right)^{i/r}\delta\left(\frac{z_2+z_0}{z_1}\right)=z_2^{-1}\left(\frac{z_1-z_0}{z_2}\right)^{-i/r}
\delta\left(\frac{z_1-z_0}{z_2}\right).$$
Now multiplying (\ref{jac}) by $z_1^{i/r}$ and taking $\Res_{z_1}$ we see
that
$$Y(Y(u,z_0)v,z_2)=\Res_{z_1}z_1^{-i/r}z^{-1}_0\delta\left(\frac{z_1-z_2}{z_0}\right)
Y(u,z_1)Y(v,z_2)$$
\begin{equation}\label{4.9}
-\Res_{z_1}z_1^{i/r}z^{-1}_0\delta\left(\frac{z_2-z_1}{-z_0}\right)
Y(v,z_2)Y(u,z_1).
\end{equation}
If $g=1$ noting that $u_nv$ is the coefficient of $z_0^{-n-1}$ in
$Y(u,z_0)v,$ multiplying (\ref{4.9}) by $z_0^{n}$ and
taking Res$_{z_0}$ we find that
\begin{eqnarray}
& &Y(u_nv,z_2)\nonumber\\
& &\ \ \ \ ={\rm Res}_{z_0}\left\{z_0^{n-1}{\rm
Res}_{z_1}\left\{\delta\left(\frac{z_1-z_2}{z_0}\right)
Y(u,z_1)Y(v,z_2)\right.\right.\nonumber\\
& &\ \ \ \ \ \ \ \ \  \left.\left.-\delta\left
(\frac{z_2-z_1}{-z_0}\right)Y(v,z_2)Y(u,z_1)\right\}\right\}\nonumber \\
& &\ \ \ \ ={\rm Res}_{z_1}\left\{(z_1-z_2)^nY(u,z_1)Y(v,z_2)-
(-z_2+z_1)^nY(v,z_2)Y(u,z_1)\right\}.\label{4.10}
\end{eqnarray}
In fact, this suggests a way to construct vertex operator algebras and their
modules for finitely generated vertex operator algebras. The vertex operator
algebras associated with highest weight representations for affine Lie algebras
and Virasoro algebra were constructed in this spirit [FZ]. For a
different construction by using so called ``the local systems''
see [L1-L2].
\end{rem}

Here we give an easy consequence of (\ref{4.10}).
In (\ref{4.10}), if $n=-1,$ we have
$$Y(u_{-1}v,z_2)=\sum_{m<0}u_mz_2^{-m-1}Y(v,z_2)+Y(v,z_2)\sum_{m\ge
0}u_mz_2^{-m-1}.$$
This suggests defining a ``normal ordering'' operation by: For $u,v\in
V,$ $m,n\in\Z,$
\begin{equation}\label{13.high54}
\mbox{$\times\atop \times$}u_mv_n\mbox{$\times\atop\times$}
=\left\{\begin{array}{ll}
u_mv_n & {\rm if}\ \ m<0\\
v_nu_m&{\rm if}\ \ m\ge 0.
\end{array}\right.
\end{equation}
Then (since $v=v_{-1}{\bf 1}$)
\begin{equation}\label{13.high55}
\ \ \ \ \ Y(u_{-1}v_{-1}{\bf 1},z_2)\!=\!Y(u_{-1}v,z_2)\!=\!\mbox{$\times\atop
\times$}Y(u,z_2)Y(v,z_2)\mbox{$\times\atop\times$}.
\end{equation}
However, this normal ordering is not in general a commutative
operation, since if $m$ and $n$ are both negative or both nonnegative,
$\mbox{$\times\atop \times$}u_mv_n\mbox{$\times\atop\times$}$ and
$\mbox{$\times\atop \times$}v_nu_m\mbox{$\times\atop\times$}$ differ
by $\pm[u_m,v_n].$ The following proposition explains that the
nilpotent property of vertex operators which holds in the algebra
will also hold for modules under a mild assumption [DL3]. In particular,
this is true for the vertex operator algebras associated with the integrable
highest weight representations of affine Lie algebras.
\begin{prop}\label{13.p16} Let $V$ be any vertex operator algebra and let
$(W,Y_W)$ be any $V$-module. Let $v\in V$ be such that the component
operators $v_n$ ($n\in{\Z}$) of $Y_W(v,z)$ all commute with one another,
so that $Y_W(v,z)^N$ is well defined on $W$ for $N\in{\N}.$ Then
\begin{equation}\label{13.13.70}
Y_W((v_{-1})^N{\bf 1},z)=Y_W(v,z)^N.
\end{equation}
In particular, if $(v_{-1})^N{\bf 1}=0$ for a fixed $N,$ then
\begin{equation}\label{13.13.71}
Y_W(v,z)^N=0.\ \ \ \ \ \Box
\end{equation}
\end{prop}

Let $V$ be a vertex operator algebra and $g$ a finite order automorphism. We
say that $V$ is
$g$-rational if $V$ has only finitely many irreducible $g$-twisted
modules and if any $g$-twisted $V$-module is completely reducible.
$V$ is rational if $V$ is $id$-rational. $V$ is holomorphic if $V$ is
rational and $V$ is the only irreducible module for itself.

\section{Examples}
\setcounter{equation}{0}
In this section we will present some well-known examples of vertex operator
algebras and their modules. We refer the reader to [D2], [DL4], [FFR],
[FLM1-FLM3], [Le], [L2], [DM2-DM3]  for various examples of twisted modules.

{\bf Vertex operator algebras associated with Heisenberg algebras.}
Let ${\bf h}$ be a vector space equipped with a symmetric nondegenerate
bilinear form $\<\cdot,\cdot\>.$ So we can identify ${\bf h}$ with its dual
${\bf h}^*$ naturally.

Viewing ${\bf h}$ as an abelian Lie
algebra, we consider the corresponding affine Lie algebra
\medskip
\begin{equation}\label{hei1}
\hat{\bf h}={\bf h}\otimes{\C}[t,t^{-1}]\oplus
{\C}c
\end{equation}
with structure defined by
\begin{equation}
\ \ \ [x\otimes t^m, y\otimes t^n]=\langle x,y\rangle m\delta_{m+n,0}c\ \
\mbox{for} \  x,y\in{\bf h},\  m,n\in{{\Z}},
\end{equation}
\begin{equation}
[c,\hat{\bf h}]=0.
\end{equation}
Set
\medskip
\begin{equation}\label{hei2}
\hat{\bf h}^+={\bf h}\otimes t{\C}[t],\ \  \hat{\bf h}^-={\bf
h}\otimes t^{-1}{\C}[t^{-1}].
\end{equation}
\medskip
The subalgebra
\begin{equation}
\hat{\bf h}_{\Z}=
\hat{\bf h}^+\oplus\hat{\bf h}^-\oplus {\C}c
\end{equation}
of $\hat{\bf h}$ is a Heisenberg algebra, in the sense that its
commutator subalgebra
coincides with its center, which is one-dimensional.
Let $\l\in {\bf h}$
and consider the
induced $\hat {\bf h}$-module
\begin{equation}\label{h1}
M(1,\l)=
U(\hat{\bf h})\otimes_{U({\bf h}\otimes{{\C}}[t]
\oplus{{\C}}c)}{{\C}}_{\l}\simeq S(\hat{\bf h}^-)\ \ \ (\mbox{linearly}),
\end{equation}
${\bf h}\otimes t{{\C}}[t]$ acting trivially on ${{\C}},$
${\bf h}$ acting as $\<h,\l\>$ for $h\in{\bf h}$ and $c$
acting as 1. We shall write $M(1)$ for $M(1,0).$
For $\alpha\in{\bf h}$ and $n\in {\Z}$ we write $\alpha(n)$
for the operator $\alpha\otimes t^n$ and  set
$$\alpha(z)=\sum_{n\in{{\Z}}}\alpha(n)z^{-n-1}$$
and define
$$
Y(v,z)=\mbox{$\circ\atop\circ$}\left(\frac{1}{(n_1-1)!}\left(\frac{
d}{dz}\right)^{n_1-1}\alpha_1(z)\right)\cdot\cdot\cdot\left(\frac{1}{(n_k-1)!}\left(\frac{d}{dz}\right)^{n_k-1}\alpha_k(z)\right)\mbox{$\circ\atop\circ$},$$
where
$$v=\alpha_1(-n_1)\cdot
\cdot\cdot\alpha_k(-n_k)\in M(1) $$
for $\alpha_1, ..., \alpha_k \in
{\bf h},\  n_1, ..., n_k \in{\Z}\ \ (n_i>0)$ and where we use a normal
ordering procedure, indicated by open colons, which signify
that the expression is to be reordered if necessary so that
all the operators
$\alpha(n)$ $(\alpha\in{\bf h},\  n<0),$ are to be placed to
the left of all the operators $\alpha(n)$
$(\alpha\in {\bf h},\ n\ge 0)$ before the expression is evaluated.
 Extend this definition to all
$v\in M(1)$ by linearity. Set ${\bf 1}=1,$
$\omega=\frac{1}{2}\sum_{i=1}^d\alpha_i(-1)^2\in M(1)$ where $\{\alpha_i\}$
is an orthonomal basis of ${\bf h}.$ The following result can be found
in [Gu].
\begin{prop} The space $M(1)=(M(1),Y,{\bf 1},\o)$ is a vertex operator
algebra and $M(1,\l)=(M(1,\l),Y)$ is a complete list
of inequivalent irreducible module for $M(1)$
for $\l\in{\bf h}.$
\end{prop}

\begin{rem} It is easy to see that the vertex operator algebra $M(1)$ is
generated by
$S=\{\alpha_i(-1)\}.$ If we identify $M(1)$ with a symmetric algebra
$\C[x_{i,n}|i=1,...,d, n=1,2,...].$ Then it is easy to see that all the vertex
operators
are built from the operators $x_{i,n}$ which is a multiplication operator
on $M(1)$ and $\frac{\partial}{\partial x_{i,n}}.$
\end{rem}

{\bf Vertex operator algebras associated with even positive definite lattices.}
Let $L$  be an even lattice, i.e., a finite-rank free abelian
group equipped  with a symmetric nondegenerate
${\Q}$-valued ${\Z}$-bilinear form\ $ \langle \cdot,\cdot\rangle ,$ not
necessarily positive definite
such that
$\langle\alpha,\alpha\rangle\in 2\Z$ for all $\alpha\in L.$
Set vector space ${\bf h}={\C}\otimes_{\Z}L$ and extend the form
$\langle\cdot,\cdot
\rangle $ from $L$ to ${\bf h}$ by ${\C}$-linearity.

The dual lattice $L^{\circ}$ of $L$ is defined to be
\begin{equation}\label{dual}
L^{\circ}=\{\beta\in {\bf h} |\<\beta, L\>\subset \Z\}.
\end{equation}
Then $L^{\circ}$ is a rational lattice  whose rank is equal to the rank of
$L.$
Let $\hat{L}^{\circ}$ be a central extension of $L^{\circ}$:
\begin{equation}
1\ \rightarrow \ \langle \omega_q\rangle
\ \rightarrow \ \hat{L}^{\circ} \ \bar{\rightarrow} \ L^{\circ} \ \rightarrow \
1,
\end{equation}
with the commutator map $c(\alpha,\beta)$ for $\alpha,\beta\in L^{\circ}$
such that
$c(\alpha,\beta)=(-1)^{\langle\alpha,\beta\rangle}$
if $\alpha,\beta\in L$,
where $\langle\omega_q\rangle $ is the cyclic group generated by a primitive
$q^{th}$
root of unity $\omega_q\in{\C}^{\times}$ and $q$ is a positive even integer.

Form the induced $\hat{L}^{\circ}$-module and ${\C}$-algebra
\begin{equation}
\begin{array}{l}
{{\C}}\{L^{\circ}\} ={\C}[\hat{L}^{\circ}]\otimes_{{\C}[\langle \omega_q\rangle
]}
{\C}\simeq{{\C}}[L^{\circ}]\ \ (\mbox{linearly}),
\end{array}
\end{equation}
where ${{\C}}[\cdot]$ denotes group algebra,
$\omega_q$ acts on ${\C}$ as multiplication by  $\omega_q$.
For $a\in\hat{L}^{\circ}$,
write $\iota(a)$ for the image of $a$ in ${{\C}}\{L^{\circ}\}$. Then the action
of
$\hat{L}^{\circ}$ on ${{\C}}\{L^{\circ}\}$ is given by:
\medskip
\begin{equation}
a\cdot\iota(b)=\iota(ab),
\end{equation}
\begin{equation}
\omega_q\cdot\iota(b)=\omega_q\iota(b)
\end{equation}
for $a,b\in\hat{L}^{\circ}.$ We give ${\C}\{L^{\circ}\}$ the ${\C}$-gradation
determined by:
\begin{equation}\label{wt}
\mbox{wt}\,\iota(a)=\frac{1}{2}\langle \bar{a},\bar{a}\rangle  \ \ \ \
\mbox{for}\ \
a\in \hat{L^{\circ}}.
\end{equation}
Also define an action of ${\bf h}$ on ${\C}\{L^{\circ}\}$ by:
\begin{equation}
h\cdot\iota(a)=
\langle h,\bar{a}\rangle \iota(a)
\end{equation}
for $h\in{\bf h}$ .
Define
\begin{equation}
z^h\cdot\iota(a)=z^{\langle h,\bar{a}\rangle }\iota(a)
\end{equation}
for $h\in{\bf h}.$
Set
\begin{equation}
V_{L^{\circ}}=M(1)\otimes_{\C}{\C}\{L^{\circ}\}\simeq S(\hat{\bf h}^-)\otimes
{\C}[L^{\circ}] \ \ \ \ (\mbox{linearly})
\end{equation}
Then $\hat{L}^{\circ},$
$\hat{\bf h}_{\Z},$ $h,$ $z^h$ $(h\in{\bf h})$ act naturally on
$V_{L^{\circ}}$ by acting on either $M(1)$ or
${{\C}}\{L^{\circ}\}$ as indicated above. In particular, $c$ acts as 1.

For a subset $M$ of $L^{\circ}$ (not necessarily a sublattice), we write
\begin{equation}\label{a3.14}
\hat M=\{b\in \hat{L}^{\circ}\ |\ \bar{b}\in M\}
\end{equation}
and
\begin{equation}\label{3.15}
{\C}\{M\}=\mbox{span}\{\iota(b)\ |\ b\in \hat M\}\subset {\bf
C}\{L^{\circ}\},
\end{equation}
\begin{equation}\label{3.16}
V_M=M(1)\otimes{\C}\{M\}\subset V_{L^{\circ}}.
\end{equation}
(Later we shall take $M$ to be a coset of a lattice.) We may apply the
considerations of the present section to $\hat{L}$ in
place of $\hat{L}^{\circ}.$ In particular, we have
\begin{equation}
{\C}\{L\}={\C}[\hat{L}]\otimes_{{\C}[\omega_p]}{\C}
\end{equation}
and
\begin{equation}
V_{L}=M(1)\otimes{\C}\{L\}.
\end{equation}

Let
\begin{equation}
L^{\circ}=\bigcup_{i\in L^{\circ}/L}(L+\lambda_i)
\end{equation}
be the coset decomposition such that $\lambda_0=0.$ Let $\Lambda_i\in
\hat L^{\circ}$ so that $\bar \Lambda_i=\lambda_i.$ Then
\begin{equation}
\hat L^{\circ}=\bigcup_{i\in L^{\circ}/L}\hat L\Lambda_i
\end{equation}
is the coset decomposition of $\hat{L}^{\circ}.$ For $i\in L^{\circ}/L,$ set
\begin{equation}
V(i)=V_{L+\lambda_i}=M(1)\otimes{\C}\{L+\lambda_i\}\simeq S(\hat{\bf
h}^-)\otimes {\C}[L+\lambda_i]\ \ \ \mbox{(linearly)}.
\end{equation}
Then we have
\begin{equation}
V_{L^{\circ}}=\coprod_{i\in L^{\circ}/L}V(i).
\end{equation}

We shall next define the untwisted vertex operators $Y(v,z)$ for $v\in
V_{L^{\circ}}.$ For $\alpha\in{\bf h}$ and $n\in {\Z}$ we write $\alpha(n)$
for the operator $\alpha\otimes t^n$ and  set
\begin{equation}\label{defp}
\alpha(z)=\sum_{n\in{{\Z}}}\alpha(n)z^{-n-1}.
\end{equation}
As before we use a normal order notation
\mbox{$\circ\atop\circ$}$\cdot$\mbox{$\circ\atop\circ$}
to signify
that the enclosed expression is to be reordered if necessary so that
all the operators
$\alpha(n)$ $(\alpha\in{\bf h},\  n<0),$ $a\in \hat{L}^{\circ}$ are to be
placed to
the left of all the operators $\alpha(n),$ $z^{\alpha}$ $
(\alpha\in {\bf h},\ n\ge 0)$ before the expression is evaluated.
For $a \in \hat{L}^{\circ}$, set
\begin{equation}\label{def1}
Y(\iota(a),z)=Y(a,z)=\mbox{$\circ\atop\circ$}e^{\int(\bar{a}(z)-
\bar{a}(0)z^{-1})}az^{\bar{a}}\mbox{$\circ\atop\circ$},
\end{equation}
using an obvious formal integration notation. Let
$$a \in \hat{L},\ \alpha_1, ..., \alpha_k \in
{\bf h},\  n_1, ..., n_k \in{\Z}\ \ (n_i>0)$$
and set
\begin{equation}\label{span}
\begin{array}{l}
v=\alpha_1(-n_1)\cdot
\cdot\cdot\alpha_k(-n_k)\otimes\iota(a) \\
\ \,\ =\alpha_1(-n_1)\cdot
\cdot\cdot\alpha_k(-n_k)\cdot\iota(a)\in V_{L^{\circ}}.
\end{array}
\end{equation}

We define
\begin{equation}\label{def2}
Y(v,z)=\mbox{$\circ\atop\circ$}\left(\frac{1}{(n_1-1)!}\left(\frac{
d}{dz}\right)^{n_1-1}\alpha_1(z)\right)\cdot\cdot\cdot\left(\frac{1}{(n_k-1)!}\left(\frac{d}{dz}\right)^{n_k-1}\alpha_k(z)\right)Y(a,z)\mbox{$\circ\atop\circ$}.
\end{equation}
This gives us a well-defined linear map
\begin{equation}
\begin{array}{l}
V_{L^{\circ}}\rightarrow(\mbox{End}\,V_{L^{\circ}})\{z\} \\
\ v\mapsto \displaystyle{Y(v,z)=\sum_{n\in{\C}}v_nz^{-n-1}, \ \
(v_n\in\mbox{End}\,V_{L^{\circ}}),}
\end{array}
\end{equation}where for any vector space $W$, we define $W\{z\}$ to be the
vector
space of $W$-valued formal series in
$z$, with arbitrary complex powers of $z$ allowed:
\begin{equation}
W\{z\}=\left\{\sum_{n\in {\C}}w_nz^n\ |\ w_n\in W\right\}.
\end{equation}
We call  $Y(v,z)$ the $untwisted$ $vertex$ $operator$ $associated$ $with$
$v.$

\begin{th} 1. $(V(0),Y,{\bf 1},\o)$ is a simple vertex operator algebra
(see [B1],[FLM3]).

2. $V(i)$ for $i\in L^{\circ}/L$ is a complete list of inequivalent
irreducible modules for $V(0)$ (see [D1], [FLM3]).

3. Any $V(0)$-module is completely reducible, that is, $V(0)$ is
rational (see [Gu]).
\end{th}

Note that from this theorem, $V_L$ is holomorphic if and only if
$L$ is self dual: $L=L^{\circ}.$ So a holomorphic vertex operator
algebra is also called a self dual vertex operator algebra [Go].

{\bf Vertex operator algebras associated with the highest weight
representations of
affine Lie algebras.} Let $\frak g$ be a simple Lie algebra over $\C,$
$\frak h$ its Cartan subalgebra and $\Delta$ the corresponding root
system. We fix a set of positive root $\Delta_+$ and a nondegenerate
symmetric invariant bilinear form $(\cdot,\cdot)$
of $\frak g$ such that the square length
of a long root is 2. The affine Lie algebra $\hat \frak g$ is defined as
$$\hat \frak g=g\otimes \C[t,t^{-1}]\oplus\C c$$
with Lie algebra structure
$$[x\otimes t^m, y\otimes t^n]=[x,y]\otimes t^{m+n}+(x,y)m\delta_{m+n,0}c$$
for $x,y\in\frak h$ and $m,n\in\Z$ where $c$ is an center element. For
convenience we shall write $x(n)$ for $x\times t^n.$
Then $\hat{\frak g}_{\pm}=\hat{\frak g}\otimes t^{\pm}{\C}[t^{\pm}]$ and
$\frak g$ (identified with with $\frak g\otimes 1$) are subalgebras.

Let $V$ be a $\frak g$-module which is extended to a $\hat\frak g_+\oplus \frak
g\oplus \C c$-module by letting $\hat \frak g_+$ act as 0 and $c$ act as
a fixed scalar $k\in \C.$ Let $\hat V_k=U(\hat\frak g)\otimes_{U(\hat\frak
g_+\oplus \frak g\oplus \C c)}V$ be the induced $\hat\frak g$-module of level
$k.$ If $V=V(\lambda)
$ is an irreducible highest weight module for $\frak g$ with
highest weight $\lambda\in\frak h^*$ we denote the corresponding $\hat V_k$ by
$\hat V_{k,\lambda}.$ Then $\hat V_{k,\lambda}$ is a Verma module. Note that
$\hat V_k=\hat V_{k,0}$ is linearly isomorphic to $U(\hat\frak g_{-})$
as vector
spaces. Our goal next is to define a vertex operator algebra structure
on $\hat V_k$ if $k$ is not equal to the dual Coxeter number and also define an
action of $\hat V_k$ on $\hat V_{k,\lambda}.$

First for $u(-1)\in \hat V_k$ we define
\begin{equation}\label{affine1}
Y(u(-1),z)=u(z)\sum_{n\in \Z}u(n)z^{-n-1}.
\end{equation}
Then using the $L(-1)$-derivation property (\ref{3.41}) to define
\begin{equation}\label{affine2}
Y(u(-n-1),z)=\frac{1}{n!}\frac{d^n}{dz^n}Y(u(-1),z).
\end{equation}
In general if $Y(v,z)$ has been defined, we use (\ref{4.10})
define $Y(u(-n)v,z)$ for $u\in \frak g$ and $0<n\in \Z$ as
\begin{equation}\label{affine3}
Y(u(-n)v,z_2)={\rm Res}_{z_1}\left\{(z_1-z_2)^{-n}Y(u,z_1)Y(v,z_2)-
(-z_2+z_1)^{-n}Y(v,z_2)Y(u,z_1)\right\}.
\end{equation}
Then we get a linear map $Y$ from $\hat V_k$ to $(\End \hat V_k)[[z,z^{-1}]].$
Set ${\bf 1}=1\times 1$ and $\omega=\frac{1}{2(k+h^{\vee})}\sum_{i}v_i(-1)^2$
where $h^{\vee}$ is the dual Coxeter number and
$\{v_i\}$ is an orthonomal basis of $\frak g$ with respect the
form $(\cdot,\cdot).$ The following theorem can be found in [FZ]
(also see [L1]):
\begin{th} (1) If $k\ne -h^{\vee}$ then $(\hat V_k,Y,{\bf 1},\omega)$ is
a vertex operator algebra.

(2) Define an action of $\hat V_k$ on $\hat V_{k,\lambda}$ by the the same
formulas. Then $(\hat V_{k,\lambda},Y)$ is a module for $\hat V_k.$
\end{th}

Now we turn our attention to the irreducible quotients. It is well-known that
$\hat V_{k,\lambda}$ has a unique maximal submodule $I(k,\lambda)$ whose
intersection
with $V(\lambda)$ (which can be identified with $1\otimes V(\lambda)$) is
0. Let $L(k,\lambda)$ be the corresponding irreducible quotient. Note that
if $k\ne 0$ and $k\ne h^{\vee}$ $I(k,0)$ is an ideal and $L(k,0)$ is a
quotient vertex operator algebra.
The following theorem can be also found in [FZ] (see [DL3] for a different
approach):

\begin{th} If $k$ is a positive integer, $L(k,0)$ is rational vertex operator
algebra whose irreducible modules are
$\{L(k,\lambda)|(\lambda,\theta)\leq k\}$ where $\theta$ is the
longest positive root in $\Delta.$
\end{th}

{\bf Vertex operator algebras associated with Virasoro algebras.} Recall
that the Virasoro algebra Vir has a basis $\{L_n|n\in\Z\}\cup \{c\}$ with
relation:
$$[L_m,L_n]=(m-n)L_{m+n}+\frac{m-m^3}{12}\delta_{m+n,0}c$$
and $c$ is in the center. Define two subalgebras
$${\rm Vir}^{\geq 0}=\oplus_{n=0}^{\infty}\C L_n,\ \  \
{\rm Vir}^{<0}=\oplus_{n=1}^{\infty}\C L_{-n}.$$
Given a pair of complex numbers $(k,h),$ the Verma module $V(k,h)$ is an
induced module
$$V(k,h)=U(\Vir)\otimes _{U(\Vir^{\geq 0})}\C_{k,h}\simeq U(\Vir^{<0})$$
where $\C_{k,h}=\C$ is a module for $\Vir^{\geq 0}$ such that $L_n1=0$
for $n>0$ and $L_01=h,\ c1=k.$

As in the case of affine algebra, we expect $V(k,0)$ is a vertex operator
algebra. But this is not quite true. The relations
$Y(L_{-1}{\bf 1},z)=\frac{d}{dz}Y({\bf 1},z)=0$ and
$\lim_{z\mapsto 0}Y(v,z){\bf 1}=v$ force $L_{-1}{\bf 1}=0.$ So it is natural
to consider the quotient module
$V_k=V(k,0)/U(\Vir)L_{-1}1$ instead of $V(k,0)$
where $1=1\otimes 1.$ Set ${\bf 1}=1+U(\Vir)L_{-1}1.$
Then $V_k$ has a basis
$$\{L{-n_1}\cdots L_{-n_s}{\bf 1}|n_i\in \N, n_1\geq n_2\geq\cdots \geq
n_2\geq2, s\geq 0\}.$$

Define operator
\begin{equation}\label{v1}
Y(L_{-2}{\bf 1},z)=L(z)=\sum_{n\in \Z}L_{n}z^{-n-2}.
\end{equation}
Again using the $L(-1)$-derivation property (\ref{3.41}) to define
\begin{equation}\label{v22}
Y(L_{-n-2}{\bf 1},z)=\frac{1}{n!}\frac{d^n}{dz^n}L(z).
\end{equation}
and use (\ref{4.10})
\begin{equation}\label{v3}
Y(L(-n)v,z_2)={\rm Res}_{z_1}\left\{(z_1-z_2)^{-n+1}L(z_1)Y(v,z_2)-
(-z_2+z_1)^{-n+1}Y(v,z_2)L(z_1)\right\}
\end{equation}
if $Y(v,z)$ has been defined. This yields
a linear map $Y$ from $V_k$ to $(\End V_k))[[z,z^{-1}]].$
Set $\omega=L_{-2}{\bf 1}.$

\begin{th} (1) $(V_k,Y,{\bf 1},\omega)$ is vertex operator algebra.

(2) Define an action of $V_k$ on $V(k,h)$ by the the same formulas.
Then $(V(k,h),Y)$ is a module for $V_k.$
\end{th}

This theorem was proved in [FZ]. Also see [H1].
Now it is easy corollary of the general theory on local systems [L1].

Now we denote by $W(k,h)$ the irreducible quotient of $V(k,h)$ module
the unique maximal submodule whose intersection with $C_{k,h}$ is 0.
Then as before $W(k,0)$ is a vertex operator algebra. Next we shall
discuss the rationality of $W(k,0).$ This is related to the discrete series
of unitary representations of the Virasoro algebra ([FQS] and [GKO]).

The work in [FQS] and [GKO] gives a
complete classification of unitary highest
weight representations of the Virasoro algebra. The highest weight
representation $W(k,h)$ is unitary
 if and only if either $(k,h)$ satisfies $k\ge 1$ and $h\ge 0,$ or
else  $(k,h)$ is among the following list:
$$k=c_m =1-\frac {6}{(m+2)(m+3)} \ \ \ (m=0,1,2\cdots ),$$
$$h=h_{r,s}^{m}=
 \frac {[(m+3)r-(m+2)s]^2 -1 } {4(m+2)(m+3)} \ \ \
 (r,s \in {\N} , 1\le s \le r \le m+1 ).$$
The unitary representations $L(c_m,h_{r,s}^m)$ for $(c_m,h_{r,s}^m)$
 in the discrete series as above are
 called the discrete series of the Virasoro algebra. If $m=0$, the only unitary
representation
 is the trivial representation $L(0,0)$. If $m=1$, there are three
 unitary representations $L(\frac{1}{2} ,0)$, $L(\frac{1}{2} ,\frac{1}{2})$
 and $L(\frac{1}{2} ,\frac{1}{16}).$

\begin{th} $W(k,0)$ is a rational vertex operator algebra if and only if
$k=c_m.$ In this case all irreducible modules are given by
$L(c_m,h_{r,s}^m).$
\end{th}

The rationality of $L(c_m,0)$ was proved  and all its irreducible modules
were found in [DMZ]. This theorem was proved completely later in [W].

\section{Twisted modules and associative algebras}
\setcounter{equation}{0}
In [Z], Zhu introduced an associative algebra $A(V)$ associated to a VOA $V$
which is extremely useful in study the representation theory of $V.$
These are analogues for $g$-twisted modules in [DLM1]. We shall review
these results form [DLM1].

Fix a VOA $V$ and an automorphism $g$ of finite $r.$ For homogeneous $u\in V^j$
and $v\in V,$ define a product $u*v$ as follows:
\begin{equation}\label{5.1}
u*v=\left\{
\begin{array}{ll}
{\rm Res}_z(Y(u,z)\frac{(z+1)^{{\rm wt}\,u}}{z}v)
 =\sum_{i=0}^{\infty}{{\rm wt}\,u\choose i}u_{i-1}v  & {\rm if}\ j=0\\
0  & {\rm if}\ j>0.
\end{array}\right.
\end{equation}
Then extends (\ref{5.1}) to a linear product $*$ on $V.$

Recall from (\ref{dec}) that $V^j$ is the eigenspace of
$g$ with eigenvalue $e^{2\pi ij/r}$ where $0\leq j<r.$
Define a subspace $O_g(V)$ of $V$ to be the linear span of all elements
$u\circ_gv$ of the following type if $u$ is homogeneous and $u\in V^j$
and $v\in V,$
\begin{equation}\label{5.3}
u\circ_gv=\left\{
\begin{array}{ll}
{\rm Res}_z(Y(u,z)\frac{(z+1)^{{\rm wt}\,u}}{z^2}v) & {\rm if}\ j=0\\
{\rm Res}_z(Y(u,z)\frac{(z+1)^{\wt\,u+j/r-1}}{z}v) & {\rm if}\ j>0.
\end{array}\right.
\end{equation}
Note that if $0<j<r,$ $u\circ_g{\bf 1}=u.$ Thus we have
\begin{lem}\label{l5.0}
If $V^{\<g\>}$ is the sub VOA of $g$-invariants of $V$ then
$V=V^{\<g\>}+O_g(V).$
\end{lem}

Let $(M,Y_g)$ be a weak $g$-twisted $V$-module.
Consider for each $ \lambda \in \C $ the subspace $
M(\lambda)=\sum_{i\in \Z}M_{\lambda+\frac{i}{K}}$ which, in fact, is a
$g$-twisted submodule of $ M$ by (\ref{1/2}). Then we have $ M =
\oplus_{\lambda \in \C/\frac{\Z}{K}} M(\lambda).$ Thus it is enough to
study the twisted module of type:
$$M=\coprod_{n\in\frac{1}{K}\Z,n\geq 0}M_{c+n}$$
where $c\in \C$ is a fixed number and $M_c\neq 0.$ We will call $M_{c}$ the
``top level''
of $M.$ For homogeneous $u\in V,$
the component operator $u_{\wt\,u-1}$ preserves each homogeneous subspace
of $M$ and in particular acts on the top level $M_c$ of $M.$
Let $o_g(u)$ be the restriction of $u_{\wt\,u-1}$ to $M_c,$ so that
we have a linear map
\begin{equation}\label{5.4}
\begin{array}{lll}
V &\to & \End(M_c)\\
u&\mapsto &o_g(u).
\end{array}
\end{equation}
Note that if $u\in V^{j/r}$ with $0<j<r,$ $o_g(u)=0$ from (\ref{1/2}).
Set $A_g(V)=V/O_g(V).$ Then we have [DLM1]:
\begin{th}\label{t5.1} (i) $A_g(V)$ is an associative algebra with
multiplication $*$ and the centralizer $C(g)$ of $g$ in $\Aut(V)$ induces
a group of algebra automorphisms of $A_g(V).$

(ii) $u\mapsto o_g(u)$ gives a representation of $A_g(V)$ on $M_c.$
Moreover, if any weak $g$-twisted module is completely reducible,
$A_g(V)$ is semisimple.

(iii) $M\to M_c$ gives a bijection between the set of equivalence
classes of simple weak $g$-twisted $V$-modules and the
 set of equivalence
classes of simple $A_g(V)$-modules.
\end{th}

The associative algebra $A(V)=A_1(V)$ was discovered in [Z] and the theorem
above was also established in this case. It was proved in [FZ] that
for the vertex operator algebra $\hat V_k,$ the associative algebra
$A(\hat V_k)$ is canonically isomorphic to $U(\frak g)$ and
$A(L(k,0))$ is isomorphic to $U(\frak g)/\<e_{\theta}^k\>$ where
$0\ne e_{\theta}\in \frak g_{\theta}$ and $\<e_{\theta}\>$ is the two-sided
ideal generated by $e_{\theta}^{k+1}.$  Moreover, $A(V_k)$ is isomorphic
to $\C[x].$ The reader can verify that $A(M(1))$ is isomorphic to
$\C[x_1,...,x_d].$

\section{Some remarks}

Finally we want to comment briefly on orbifold theory and on
generalizations of vertex operator algebras.

It is believed that the module category
for a rational VOA is equivalent to the module category of a
Hopf algebra associated with the vertex operator algebra. In the
case of VOAs associated with the level $k$
highest weight unitary representations
for affine Lie algebras $\hat \frak g,$ the Hopf algebras are essentially the
quantum group $U_q(\frak g)$ where $q=e^{2\pi i/h+k}$ where $h$ is the
dual Coxeter number. If $V$ is holomorphic the corresponding Hopf algebra
is $\C.$ It is natural next to determine the Hopf algebra for
$V^G$ where $V$ is a holomorphic VOA, $G$ is a finite automorphism
group and $V^G$ is the $G$-invariants which is a vertex operator
subalgebra of $V.$ This is the so-called orbifold theory in the physical
literature. In fact, the moonshine module
is an $\Z_2$-orbifold theory constructed from the vertex operator algebra
associated with Leech lattice and an order 2 automorphism induced from
$-1$ isometry of the Leech lattice (a $\Z_p$-orbifold construction of the
moonshine module have been studied in [DM1] and [Mon]).

The main new feature of the orbifold theory is the introduction of twisted
modules because any $g$-twisted module
is an ordinary module for $V^G$ for $g\in G.$ It is conjectured that for
holomorphic vertex operator algebra $V$ and finite automorphism group
$G,$ the corresponding Hopf algebra is the quasi quantum double
$D^{c}(G)$ where $c\in H^3(G,\C)$ is uniquely determined by $V$ (see
[DVVV], [DPR], [DM2-DM5]). There have been a lot progress in
this direction in [DM4-DM5] for nilpotent group $G.$
The results concerning the modular invariance of trace functions
(or correlation functions on the torus) developed in [Z] and [DM6] play
important roles in the orbifold theory. (These results also explain
mysterious relations among affine Lie algebras, Virasoro algebra, monster
group and the modular group.)

In [DL1-DL3], the theory of vertex operator algebras has been generalized
in a systematic way at three successively more general levels, all of which
incorporate one-dimensional braid group representations intrinsically into
the algebraic structure: First, the notion of ``generalized vertex operator
algebra'' incorporates such structure as $Z$-algebras, paraferminon algebras,
and the vertex operator superalgebras. Next, what we term ``generalized
vertex algebras'' further encompass the algebras of vertex operators
associated with rational lattices. Finally, the most general of the
three notions,
that of ``abelian intertwining algebras,'' also include
intertwining operators for
certain classes of vertex operator algebras. See [L4], [H4] and
[DLM2] for more examples of abelian intertwining algebras related to the
simple currents and the moonshine module.
The notion of generalized vertex algebra has been also
independently introduced and studied in [FFR] with different motivations,
examples (involving spinor constructions) and axiomatic approach from ours,
and also see [Mos]. See also [LZ] and [H2] for the notion of topological
vertex algebra.

\noindent Address: {Department of Mathematics,
 University of California,
 Santa Cruz, CA 95064}

\hspace{0.8cm} dong@dong.ucsc.edu

\end{document}